\documentclass[reqno,12pt]{amsart}

\usepackage{color}
\usepackage{enumerate}

\numberwithin{equation}{section}

\theoremstyle{plain}
\newtheorem{thm}{Theorem}[section]

\newtheorem{hyp}[thm]{Hypotheses}

\theoremstyle{remark}

\theoremstyle{definition}

\def\question#1{\ifmmode\text{\bf \color{red} Q: #1}\else{\bf
    Q}\footnote{#1}\fi:\color{red}}

\usepackage{amscd,amssymb,comment,epic,eepic,euscript,graphics}

\newcommand{\la}{\lambda}

\newcommand\Mcal{{\mathcal{M}}} 

\newcommand\Lcal{{\mathcal{L}}}
\newcommand\Ical{{\mathcal{I}}}
\newcommand\Reals{{\mathbb R}}
\newcommand\Complex{{\mathbb C}}
\newcommand\Nats{{\mathbb N}}

\newcommand{\BH}{\mathcal{B}(\mathcal{H})}
\newcommand{\Hcal}{\mathcal{H}}

\newcommand\Tr{{\mathrm{Tr}}}

\newcommand{\im}{\text{\rm Im}}
\newcommand{\re}{\text{\rm Re}}

\begin{document}

\title{Taylor approximations of operator functions}

\author[Skripka]{Anna Skripka$^{*}$} \email{skripka@math.unm.edu}

\thanks{\indent\llap{${}^{*}$}Research supported in part by NSF
  grant DMS-1249186.}

\address{A.S., Department of Mathematics and Statistics, University of New Mexico, 400 Yale Blvd NE, MSC01 1115, Albuquerque, NM 87131, USA}




\begin{abstract}
This survey on approximations of perturbed operator functions addresses recent advances and some of the successful methods.
\end{abstract}

\maketitle

\section{Introduction}
An active mathematical investigation of perturbed operator functions started in as early as 1950's, following a series of physics papers by I.~M.~Lifshits on the change of the free energy of a crystal due to appearance of a small defect. The latter research in physics gave birth to the Lifshits-Krein spectral shift function \cite{L,Krein,KreinUn}, which has become a fundamental object in perturbation problems of mathematical physics. Subsequent attempts to include more general perturbations than those in \cite{Krein,KreinUn} have resulted in consideration of higher order Taylor approximations of perturbed operator functions and introduction of Koplienko's higher order spectral shift functions \cite{Kop,Neidhardt,ds,PSS,PS-circle,PSS-circle}.

Approximation of operator functions also arises in problems of noncommutative geometry involving spectral flow (see, e.g., \cite{sf}) and spectral action functional (see, e.g., \cite{CC}). This investigation was initially carried out independently of the study of the spectral shift functions. However, a recent unified approach to the Lifshits-Krein spectral function and the spectral flow allowed to establish that these two objects essentially coincide \cite{ACS}. Higher order Taylor formulas have been derived for spectral actions in \cite{vanS}, with restrictions on the operators relaxed in \cite{as} by applying more universal perturbation theory techniques.

The proof of existence of the first order (Lifshits-Krein) spectral shift function, which is due to M.~G.~Krein, relied on the theory of analytic functions and was of a different nature than the proofs of the other mentioned results on the approximations of operator functions. An important object in higher order Taylor approximations is the G\^{a}teaux derivative of an operator function. When the initial operator and the perturbation do not commute, the G\^{a}teaux derivative is a complex object, whose complexity increases with the order of differentiation. Treatment of such derivatives and subsequent derivation of Taylor approximations was based on a delicate noncommutative analysis, which had been developing for some $60$ years.

To proceed to a detailed discussion of the aforementioned and further results and methods, we need to fix some notation. We work with a pair of operators defined in a separable Hilbert space $\Hcal$, denoting the initial operator $H_0$ and its perturbation $V$. The perturbation is always a bounded operator and, moreover, some summability restrictions are imposed either on $V$ or $H_0$. In some instances, $H_0$ is allowed to be unbounded, and we will consider only closed densely defined unbounded operators. For sufficiently nice scalar functions $f$, we consider the operator functions $f(H_0)$ and $f(H_0+V)$ given by the functional calculus. We are interested in some scalar characteristics associated with perturbations that are calculated using traces (a canonical trace $\Tr$, a Dixmier trace $\Tr_\omega$, a normal trace on a semi-finite von Neumann algebra $\tau$, and, more generally, any trace $\tau_\mathcal{I}$ on a normed ideal $\mathcal{I}$ continuous in the ideal norm).

We consider the remainders of the Taylor approximations
\begin{align*}
R_{n,H_0,V}(f)=f(H_0+V)-\sum_{k = 0}^{n-1} \frac 1{k!}\, \frac {d^k}{dt^k}\bigg|_{t = 0} f(H_0 + tV),
\end{align*}
where $n\in\Nats$ and the G\^{a}teaux derivatives $\frac{d}{dt}\big|_{t=0}f(H_0+tV)$ are evaluated in the uniform operator topology. If the $n$-th order G\^{a}teaux derivative is continuous on $[0,1]$, then we have the integral representation
for the remainder
\begin{align}
\label{ir}
R_{n,H_0,V}(f)=\frac{1}{(n-1)!}\int_0^1 (1-t)^{n-1}\frac {d^n}{dt^n} f(H_0 + tV)\,dt,
\end{align}
which can be proved by applying functionals in the dual space $(\mathcal{B}(\Hcal))^*$ of the algebra of bounded linear operators on $\Hcal$ and reducing the problem to the scalar case. The questions we are interested in consist in establishing more specific properties of the remainders
$R_{n,H_0,V}(f)$.

\section{Schatten class perturbations}

In this section, we discuss Taylor approximations in the classical setting of perturbations belonging to the Schatten-von Neumann ideals of compact operators
\[S^\alpha=\big\{A\in\BH:\; \|A\|_n:=\big(\Tr(|A|^\alpha)\big)^\frac1\alpha<\infty\big\},\quad \alpha\in [1,\infty)\]
(see, e.g., \cite{Simon}).
The operator functions under consideration come from either polynomials $\mathcal{P}$ or the functions with nice Fourier transforms
\[\mathcal{W}_n=\big\{f:\; f^{(j)}, \widehat{f^{(j)}}\in L^1(\Reals), j=0,\ldots,n\big\}.\]
The class $\mathcal{W}_n$ includes such widely used sets of functions as $C_c^{n+1}(\Reals)$ and the rational functions in $C_0(\Reals)$, which we denote by $\mathcal{R}$.

\subsection{\bf Spectral shift functions}
As a joint finding of many investigations, we have the following representation for the Taylor remainders corresponding to self-adjoint perturbations of self-adjoint operators.

\begin{thm}
\label{saSn}
If $H_0=H_0^*$ and $V=V^*\in S^n$, $n\in\Nats$, then there exists a unique real-valued function $\eta_n=\eta_{n,H_0,V}\in L^1(\Reals)$ and a constant $c_n>0$ such that \[\|\eta_n\|_1\leq c_n \|V\|_n^n\] and
\begin{equation}
    \label{nTraceFormula}
    \Tr \left( f(H_0+V)-\sum_{k = 0}^{n-1} \frac 1{k!}\, \frac {d^k}{dt^k}\bigg|_{t = 0}
      f(H_0 + tV) \right) = \int_\Reals f^{(n)}
    (t)\, \eta_n  (t)\, dt,
  \end{equation}
for $f\in\mathcal{W}_n$.
\end{thm}
The cases $n=1$, $n=2$, and $n\geq 3$ are due to \cite{Krein} (see also \cite{summer}), \cite{Kop}, and \cite{PSS}, respectively. The formula \eqref{nTraceFormula} has been extended from $\mathcal{W}_n$ to the Besov class $B_{\infty1}^n(\Reals)$ in \cite{PellerKr}, \cite{PellerKo}, and \cite{AP}, respectively. Differentiability of operator functions in the setting most applicable to Theorem \ref{saSn} is discussed in \cite{PellerMult} and \cite{Azamov0}. The results of \cite{Azamov0,PellerMult} can also be used to justify that the trace on the left hand side of \eqref{nTraceFormula} is well defined.

The function $\eta_n$ provided by Theorem \ref{saSn} is called the $n$-th order spectral shift function associated with the pair of operators $(H_0,H_0+V)$. The name to $\eta_1$ was given by M.~G.~Krein and can be understood from I.~M.~Lifshits' formula
\[\eta_1(\la)=\Tr\big(E_{H_0}((-\infty,\la))\big)-\Tr\big(E_{H_0+V}((-\infty,\la))\big)\]
holding for $H_0$ and $V$ finite matrices, where $E_H$ denotes the spectral measure of $H$.
A number of remarkable connections of the first-order spectral shift function to other objects of mathematical physics can be found in the brief survey \cite{BP}. More detailed discussion of the first order spectral shift function can be found in \cite{BY,Simon,Yafaev} and of the second order one in \cite{GPS}.
When a perturbation $V$ is in the Hilbert-Schmidt class $S^2$, the higher order spectral shift functions can be expressed via the lower order ones \cite{ds,unbdd}. The former are more sensitive to the displacement of the spectrum under perturbation, as demonstrated in \cite{convexity,moissf}.

The question of validity of
\begin{align}
\label{nsatf}
\Tr\big(f(H_0+V)-f(H_0)\big)=\int_\Omega f'(t)\,\eta_1(t)\,dt,
\end{align}
was also investigated for non-self-adjoint operators $H_0$ and $H_0+V$. Here the set $\Omega\subset\Complex$ is determined by $H_0$ and $V$.
The trace formula \eqref{nsatf} with $\Omega=\mathbb{T}$ (the unit circle) was proved in \cite{KreinUn} for unitary operators $H_0$ and $H_0+V$ such that $V\in S^1$. The case of arbitrary bounded operators $H_0$ and $H_0+V$ differing by $V\in S^1$ is naturally harder than the case of self-adjoint operators. If $H_0$ and $H_0+V$ are contractions, then \eqref{nsatf} holds with $\Omega=\mathbb{T}$ and the function $\eta_1$ replaced by a finite complex-valued measure $\nu_1$, for every $f$ analytic on a disc centered at zero of radius $r>1$ \cite{dsD}. Attempts to get more information about the structure of the measure $\nu_1$ (for example, extract an absolutely continuous component) resulted in consideration of only selected pairs of contractions
and brought to modification of \eqref{nsatf} with passage to a more general type of integration. The relevant discussion (also for dissipative operators $H_0$ and $H_0+V$) can be found in \cite{AdN,AdP,K87,MSZ,N87,N88,R84,R94,R96}.

The higher order version of \eqref{nsatf} for pairs of bounded operators has a more plausible formulation.

\begin{thm}
\label{conSn}
Let $H_0$ and $H_0+V$ be contractions and assume that $V\in S^n$, $n\geq 2$. Then, there exists a function $\eta_n=\eta_{n,H_0,V}$ in $L^1(\mathbb{T})$ such that
\begin{align}
\label{tfcon}
 \Tr \left( f(H_0+V)-\sum_{k = 0}^{n-1} \frac 1{k!}\, \frac {d^k}{dt^k}\bigg|_{t = 0}
      f(H_0 + tV) \right)=\int_\mathbb{T} f^{(n)}(z)\,\eta_n(z)\,dz,
\end{align}
for $f\in\mathcal{P}$.
Furthermore, there exists a constant $c_n>0$ such that 
a function $\eta_n$ satisfying \eqref{tfcon} can be chosen so that
\begin{align}
\label{ssfest}
\|\eta_n\|_1\leq c_n\|V\|_n^n.
\end{align}
\end{thm}

The case $n=2$ for $H_0$ and $H_0+V$ unitaries, where the derivative is evaluated along a multipicative path of unitaries instead of $t\mapsto H_0+tV$, is due to \cite{Neidhardt} (with later extension of the class of functions $f$ in \cite{PellerKo}) and for arbitrary pairs of contractions $H_0$ and $H_0+V$ joined by the path $t\mapsto H_0+tV$ is due to \cite{PS-circle}. The case $n\geq 3$ is established in \cite{PSS-circle}. The spectral shift function $\eta_n$ satisfying Theorem \ref{conSn} is determined uniquely only up to an analytic term (that is, the equivalence class of $\eta_n$ in the quotient space $L^1(\mathbb{T})/H^1(\mathbb{T})$ is uniquely determined).  Theorem \ref{conSn} can be extended to more general functions $f$. In particular, \eqref{tfcon} with $n=2$ is established for analytic functions $f$ in \cite{dsD}, as discussed in Subsection \ref{subseqgt} for more general traces.

\subsection{\bf Proof strategy}\label{subsec22}
The proofs of Theorems \ref{saSn} and \ref{conSn} are very subtle and technically involved, so we will give only a flavor of some basic ideas. For simplicity we assume that $\|H_0\|\leq 1$, $\|H_0+V\|\leq 1$, $V\in S^n$, and $f\in\mathcal{P}$. Then our goal is the formula
\begin{align}
\label{simplified}
\Tr\big(R_{n,H_0,V}(f)\big)=\int_{\mathbb{T}}f^{(n)}(z)\nu_n(z)\,dz,
\end{align}
where $\nu_n$ is a finite measure, with total variation bounded by
\begin{align}
\label{nuest}
\|\nu_n\|\leq c_n\|V\|_n^n.
\end{align}

From the integral representation for the remainder \eqref{ir}, we derive
\[\Tr\big(R_{n,H_0,V}(f)\big)=\frac{1}{(n-1)!}\int_0^1 (1-t)^{n-1}\,\Tr\left(\frac{d^n}{ds^n}\bigg|_{s=t} f(H_0 + sV)\right)dt.\]
Thus, if we prove
\begin{align}
\label{dest}
\sup_{t\in [0,1]}\left|\frac{1}{n!}\,\Tr\left(\frac{d^n}{ds^n}\bigg|_{s=t} f(H_0 + sV)\right)\right|\leq c_n \|V\|_n^n\cdot\big\|f^{(n)}\big\|_\infty,
\end{align}
then application of the Hahn-Banach theorem and the Riesz representation theorem for the dual space of $C(\mathbb{T})$ implies existence of a measure $\nu_n$ satisfying \eqref{simplified} and \eqref{nuest}.

For $n=1$, we have
\[\Tr\left(\frac {d}{ds}\bigg|_{s=t} f(H_0 + sV)\right)=\Tr\big(f'(H_0+tV)V\big),\] which in case of $f$ a polynomial follows from the straightforward calculation of the derivative and some combinatorics. Applying the H\"{o}lder and von Neumann inequalities then implies \eqref{dest} with $n=1$ and $c_1=1$. This reasoning does not allow to establish the absolute continuity of $\nu_1$ (which was established in \cite{Krein}), but it can be generalized to apply to the higher order case. If, in addition, we take $H_0$ and $V$ to be self-adjoint, then application on the spectral theory allows to derive an explicit formula for $\nu_1$, as it was done in \cite{Birman72}.

Apart from the case of commuting $H_0$ and $V\in S^2$, we do not have the convenient equality
$\Tr\left(\frac {d^2}{ds^2}\big|_{s=t} f(H_0 + sV)\right)=\Tr\big(f''(H_0+tV)V^2\big)$. However, since
the set function
$A_1\times A_2\mapsto \Tr\big(E_{H_0+tV}(A_1)V E_{H_0+tV}(A_2)V\big)$,
where $A_1,A_2$ are Borel subsets of $\Reals$, uniquely extends to a measure on $\Reals^2$ with total variation $\|V\|_2^2$, we have
\[\Tr\left(\frac{d^2}{dt^2}\bigg|_{s=t}f(H_0+sV)\right)=
\int_{\Reals^2}(f')^{[1]}(\la_1,\la_2)\,\Tr\big(E_{H_0+tV}(d\la_1)V E_{H_0+tV}(d\la_2)V \big)\]
(see, e.g., \cite[Theorem 3.12]{moissf}), which along with the estimate for the divided difference $\|(f')^{[1]}\|_\infty\leq \|f''\|_\infty$ implies \eqref{dest} with $n=2$ and $c_2=\frac12$.

When $n\geq 3$, the set function $A_1\times\cdots\times A_n\mapsto \Tr\big(E_{H_0+tV}(A_1)V\ldots E_{H_0+tV}(A_n)V\big)$ can fail to extend to a measure of finite variation on $\Reals^n$ (see \cite[Section 4]{ds}). 
This is one of the reasons suggesting that the case $n\geq 3$ requires much more delicate (noncommutative) analysis of operator derivatives than the case $n<3$.

Pioneering estimates for norms of $n$-th order operator derivatives are attributed to Yu.~L.~Daleckii and S.~G.~Krein \cite{DK}. In \cite{DK}, $H_0=H_0^*$ and $V=V^*\in\BH$, a scalar function $f$ belongs to $C^{2n}(\Reals)$, and the estimates depend on the size of the spectrum of the operator $H_0$. Development of the Birman-Solomyak double operator integration (see, e.g., \cite{BS}) and subsequent multiple operator integration (see \cite{PellerMult} and also \cite{Azamov0}) resulted in significant improvement of the estimates for operator derivatives. It follows from~\cite{PellerMult} that for $H_0=H_0^*$ and $V=V^*\in S^n$,
\begin{equation*}
    \sup_{t\in [0,1]}\left|\Tr\left(\frac {d^n}{ds^{n}}\bigg|_{s=t}f(H_0+sV)\right)\right|\leq c_n \left\|
      f \right\|_{B_{\infty1}^n(\Reals)}\cdot \left\| V\right\|_n^n,
\end{equation*} where $f\in B_{\infty1}^n(\Reals)$;
however, the norm $\left\| f \right\|_{B_{\infty1}^n(\Reals)}$ is
greater than the norm $\|f^{(n)}\|_{L^\infty(\Reals)}$. The powerful estimates \eqref{dest} are established in the following theorems.

\begin{thm} $($\cite{PSS-circle}$)$ 
\label{mecon}
If $\|H_0\|\leq 1$, $\|H_0+V\|\leq 1$, and $n\in\Nats$,
then there exists a constant $c_n>0$ such that for every $f\in\mathcal{P}$ the following estimates hold.
\begin{enumerate}[(i)]
\item If $\beta>n$ and $V\in S^\beta$, then
\begin{equation*}
\sup_{t\in [0,1]}\left\|\frac{d^n}{ds^n}\bigg|_{s=t}f(H_0+sV)\right\|_{\frac{\beta}{n}}\leq c_n \|V\|_\beta^{n}\cdot\big\|f^{(n)}\big\|_{L^\infty(\mathbb{T})}.
\end{equation*}

\item If $V\in S^n$, then
\begin{equation*}
\sup_{t\in [0,1]}\left|\Tr\left(\frac{d^n}{ds^n}\bigg|_{s=t}f(H_0+sV)\right)\right|\leq c_n \|V\|_n^n\cdot \big\|f^{(n)}\big\|_{L^\infty(\mathbb{T})}.
\end{equation*}
\end{enumerate}
\end{thm}

\begin{thm} $($\cite{PSS}$)$ 
\label{mesa}
If $H_0=H_0^*$, $V=V^*$, and $n\in\Nats$,
then there exists a constant $c_n>0$ such that for every $f\in\mathcal{W}_n$ the following estimates hold.
\begin{enumerate}[(i)]
\item If $\beta>n$ and $V\in S^\beta$, then
\begin{equation*}
\sup_{t\in [0,1]}\left\|\frac{d^n}{ds^n}\bigg|_{s=t}f(H_0+sV)\right\|_{\frac{\beta}{n}}\leq c_n \|V\|_\beta^{n}\cdot\big\|f^{(n)}\big\|_{L^\infty(\Reals)}.
\end{equation*}

\item If $V\in S^n$, then
\begin{equation*}
\sup_{t\in [0,1]}\left|\Tr\left(\frac{d^n}{ds^n}\bigg|_{s=t}f(H_0+sV)\right)\right|\leq c_n \|V\|_n^n\cdot \big\|f^{(n)}\big\|_{L^\infty(\Reals)}.
\end{equation*}
\end{enumerate}
\end{thm}

The proofs of Theorems \ref{mecon} and \ref{mesa} (and also analogous estimates for polylinear transformations more general than operator derivatives) include a subtle synthesis of advanced techniques from harmonic, functional, complex analysis and noncommutative $L^p$ spaces as well as development of a novel approach to multiple operator integration. The principal two cases here are the ones of self-adjoints and unitaries, while the case of contractions reduces to the case of unitaries by applying the Sz.-Nagy-Foia\c{s} dilation theory \cite{SzNF}.

\subsection{\bf Operator Lipschitz functions}

Derivation of the estimates of Theorems \ref{mecon} and \ref{mesa} was preceded by resolution of Krein's conjecture on whether every Lipschitz function on $\Reals$ is operator Lipschitz. Detailed discussion of the problem, including references to partial results, can be found in \cite{Pellersurvey,PS-Acta}; here we only state the concluding result and mention some generalizations.
\begin{thm}$($\cite{PS-Acta}$)$
Let $f$ be a Lipschitz function on $\Reals$. Then, for every $\alpha\in(1,\infty)$, there is a constant $c_\alpha>0$ such that
\[\|f(B)-f(A)\|_\alpha\leq c_\alpha\|B-A\|_\alpha\cdot\|f\|_{Lip},\]
for all $A=A^*$, $B=B^*$, defined in $\Hcal$ with $B-A\in S^\alpha$.
\end{thm}
The best constant $c_\alpha\sim \frac{\alpha^2}{\alpha-1}$ is obtained in \cite{Caspers}. It is known that not every Lipschitz function is operator Lipschitz in $S^1$ and in $\BH$ (i.e., when $\alpha\in\{1,\infty\}$) \cite{F1,F2,F3}. Operator Lipschitzness of functions of normal operators and of functions of several variables is discussed in \cite{APPS,KPSS}.

\section{Some natural generalizations}

If a perturbation $V$ is not compact and no additional restriction on $H_0$ is imposed, then the canonical trace $\Tr$ of $R_{n,H_0,V}(f)$ is not defined. Depending on the problem, one can consider another trace that is defined on $R_{n,H_0,V}(f)$ for rather general $H_0$, $V$, and $f$, or impose extra restrictions on $H_0$, $f$, and/or $V$ to ensure $R_{n,H_0,V}(f)\in S^1$.

\subsection{\bf Compact resolvents and similar conditions}

Perturbations that arise in the study of differential operators are multiplications by functions defined on $\Reals^d$, which are not compact operators. In this case, the condition $V\in S^n$ gets replaced by a restriction on the resolvent of the initial operator $H_0$.

If $H_0$ equals the negative Laplacian $-\Delta$ and the operator $V$ act as multiplication by a real-valued function in $L^1(\Reals^3)\cap L^\infty(\Reals^3)$, then
\begin{align}
\label{rc1}
(H_0-zI)^{-1}-(H_0+V-zI)^{-1} \in S^1,\quad z\in\Complex\setminus\Reals
\end{align}
(see, e.g., \cite{BY}).
Due to the invariance principle for the first order spectral shift function (see, e.g., \cite{BY}), the problem for a pair of self-adjoint operators $(H_0,V)$ satisfying \eqref{rc1} reduces to the problem for a pair of unitaries with difference in $S^1$, and \eqref{nTraceFormula} with $n=1$  holds for $f\in C_c^\infty(\Reals)\cup\mathcal{R}$, as established in \cite{KreinUn}. In this case, $\eta_1$ is an element of $L^1\big(\Reals,\frac{1}{1+t^2}dt\big)$. Existence of the first order spectral shift function under more general resolvent conditions is discussed in \cite{Kop77,Yafaev07}.

If $H_0=-\Delta$ and $V$ is a multiplication by a real-valued function in $L^2(\Reals^3)\cap L^\infty(\Reals^3)$, then instead of the condition \eqref{rc1}, we have
\begin{align}
\label{rc2}
(I+H_0^2)^{-1/4}V\in S^2
\end{align}
(see, e.g., \cite{HS-compatible}). It is established in \cite{Kop} that for a pair of self-adjoint operators $(H_0,V)$ satisfying \eqref{rc2}, there exists $\eta_2\in L^1\big(\Reals,\frac{1}{1+t^2}dt\big)$ such that the trace formula \eqref{nTraceFormula} with $n=2$  holds for $f\in \mathcal{R}$. A modified trace formula is obtained in \cite{HS-compatible} for a pair $(H_0,V)$ satisfying $(I+H_0^2)^{-1/2}V\in S^2$. The proofs are based on multiple operator integration techniques developed to partly compensate for the lack of the invariance principle under the assumption \eqref{rc2}.

In perturbation problems of noncommutative geometry, typical assumptions on the operators are that the resolvent of $H_0$ is compact and $V\in\BH$. The following result is obtained in \cite{as}, relaxing assumptions on $H_0$ and $V$ made in \cite{vanS}.
\begin{thm}
\label{asexp}
Let $H_0=H_0^*$ be defined in $\mathcal{H}$ and have compact resolvent and let $V=V^*\in\BH$. Let $\{\mu_k\}_{k=1}^\infty$ be a sequence of eigenvalues of $H_0$ counting multiplicity and let $\{\psi_k\}_{k=1}^\infty$ be an orthonormal basis of the respective eigenvectors. Then, for each function $f \in C_c^{n+1}(\Reals)$, with $n\in\Nats$,
\begin{align*}
&\Tr\big(f(H_0+V)\big)-\Tr\big(f(H_0)\big)\\
\nonumber
&\quad=\sum_{p=1}^{n-1} \frac1p\sum_{i_1,\ldots,i_p}(f')^{[p-1]}(\mu_{i_1},\ldots,\mu_{i_p})\,
\left<V\psi_{i_1},\psi_{i_2}\right>\cdots\left<V\psi_{i_p},\psi_{i_1}\right>+\Tr\big(R_{H_0,f,n}(V)\big),
\end{align*}
where
\[\Tr\big(R_{H_0,f,n}(V)\big)=\mathcal{O}\big(\|V\|^n\big).\]
\end{thm}
Moreover, the trace formula \eqref{nTraceFormula} with $f\in C_c^3(\Reals)$ is established in \cite{ACS} for $n=1$ (this is also discussed in the next subsection) and, under the additional assumption $(I+H_0^2)^{-1/2}\in S^2$, in \cite{as} for $n=2$. The respective spectral shift functions $\eta_1$ and $\eta_2$ are locally integrable.

Taylor asymptotic expansions and spectral distributions have also been considered in the study of pseudodifferential operators (see, e.g., \cite{BoucletKo}).

\subsection{\bf Operators in a semifinite von Neumann algebra}

Let $\Mcal$ be a semifinite von Neumann algebra of bounded linear operators defined on $\Hcal$ and let $\tau$ be a semifinite normal faithful trace $\tau$ on $\Mcal$. (The definitions can be found in, e.g., \cite{LSZ}.) Note that $(\BH,\Tr)$ is one of examples of $(\Mcal,\tau)$. Let $H_0$ be either an element of $\Mcal$ or an unbounded closed densely defined self-adjoint operator affiliated with $\Mcal$ (that is, all the spectral projections of $H_0$ are elements of $\Mcal$). The perturbation $V$ is taken to be a bounded element of the noncommutative $L^p$-space associated with $(\Mcal,\tau)$, that is,
\[V\in \Lcal^n=\big\{A\in\Mcal:\; \|A\|_n:=\tau(|A|^n)^\frac1n<\infty\big\},\quad n\in\Nats.\]

\begin{thm}
\label{saLn}
If $H_0=H_0^*$ is affiliated with $\Mcal$ and $V=V^*\in \Lcal^n$, $n\in\Nats$, then there exists a unique real-valued function $\eta_n=\eta_{n,H_0,V}\in L^1(\Reals)$ and a constant $c_n>0$ such that \[\|\eta_n\|_1\leq c_n \|V\|_n^n\] and
\begin{equation}
    \label{nTraceFormulaL}
    \tau \left( f(H_0+V)-\sum_{k = 0}^{n-1} \frac 1{k!}\, \frac {d^k}{dt^k}\bigg|_{t = 0}
      f(H_0 + tV) \right) = \int_\Reals f^{(n)}
    (t)\, \eta_n  (t)\, dt,
  \end{equation}
for $f\in\mathcal{W}_n$.
\end{thm}
The case $n=1$ was established first for a bounded operator $H_0$ in \cite{Carey} and then for an unbounded operator in \cite{Azamov}. The case $n=2$ is due to \cite{ds,unbdd} and $n\geq 3$ is due to \cite{PSS}. The strategy of the proof is as described in Subsection \ref{subsec22}; this strategy can be implemented because noncommutative $L^p$-spaces have much in common with Schatten ideals (see, e.g., \cite{PX}).

The first order spectral shift function for a pair of $\tau$-Fredholm operators differing by a $\tau$-compact perturbation is known to coincide with the spectral flow \cite[Theorem 3.18]{ACS}. It is also established in \cite{ACS} that \eqref{nTraceFormulaL} with $n=1$ holds for $H_0$ having $\tau$-compact resolvent. (In the case $(\Mcal,\tau)=(\BH,\Tr)$, a $\tau$-compact operator is merely a compact operator.)

\begin{thm}$($\cite{ACS}$)$\label{R1}
If $H_0=H_0^*$ is affiliated with $\Mcal$ and has a $\tau$-compact resolvent and if $V=V^*\in\Mcal$, then, for $f\in C_c^3((a,b))$,
\begin{align*}
\tau\big(f(H_0+V)\big)=\tau\big(f(H_0)\big)+\int_\Reals f'(\la)\tau\big(E_{H_0}((a,\la])-E_{H_0+V}((a,\la])\big)\,d\la.
\end{align*}
\end{thm}

Analogs of \eqref{nTraceFormulaL} with $n=1$ and $n=2$ for pairs of arbitrary (non-self-adjoint) operators in $\Mcal$ differing by a perturbation $V\in\Lcal^n$ are obtained in \cite{dsD}. As to the case $n\geq 3$, the results of Theorem \ref{mesa} can be extended to pairs of operators in $\Mcal$ by applying dilation of contractions in $\Mcal$ to unitary operators in semi-finite von Neumann algebras constructed in \cite{dsD}.

\subsection{\bf General traces}
\label{subseqgt}

The canonical trace $\Tr$ is widely used, but it is not the most ``typical" trace. The distinctive feature of $\Tr$ is that it is normal, i.e, has the property of monotonicity. A continuous trace on a normed ideal of compact operators in $\BH$ other than $S^1$ has a singular component, which vanishes on finite rank operators. Detailed discussion of traces and applications of singular traces to classical and noncommutative geometry can be found in \cite{LSZ}.

Let $\Mcal$ be a semifinite (von Neumann) factor and $\Ical$ a symmetrically normed ideal of $\Mcal$ with norm $\|\cdot\|_\Ical$. (The definitions can be found, e.g., in \cite{dsD,LSZ}.) Let $\tau_\Ical$ be a trace on $\Ical$ bounded with respect to the ideal norm $\|\cdot\|_\Ical$. Examples of $(\Ical,\tau_\Ical)$ include $(S^1,\Tr)$, $(\Lcal^1,\tau)$, where $\tau$ is the normal faithful semifinite trace on $\Mcal$, and $(\Lcal^{(1,\infty)},\Tr_\omega)$, where $\Lcal^{(1,\infty)}$ denotes the dual Macaev ideal and $\Tr_\omega$ the Dixmier trace on it corresponding to a generalized limit $\omega$ on $\ell^\infty(\Nats)$.

The following results are obtained in \cite{dsD}.

\begin{hyp}
\label{cases}
Consider a set $\Omega$, a closed, densely defined operator $H_0$ affiliated to $\Mcal$, an operator $V\in\Ical$
and a space $\mathcal{F}$ of functions that satisfy one of the following assertions.
\begin{enumerate}[(i)]
\item \label{i} $\Omega=\text{\rm conv}\big(\sigma(H_0)\cup\sigma(H_0+V)\big)$, $H_0=H_0^*\in\Mcal$, $V=V^*$,
$\mathcal{F}=C^3(\Reals)$;
\item \label{ii}$\Omega=\Reals$, $\im(H_0)\geq 0$, $\im(H_0+V)\geq 0$, and
\[
\mathcal{F}=\text{\rm span}\left\{\la\mapsto (z-\la)^{-k}:\, k\in\Nats,\, \im(z)<0\right\};
\]
\item \label{iii} $\Omega=\mathbb{T}$, $\|H_0\|\leq 1$, $\|H_0+V\|\leq 1$, and $\mathcal{F}$ is the set of all functions that are analytic on discs centered at $0$ and of radius strictly larger than $1$.
\end{enumerate}
\end{hyp}

\begin{thm}
\label{mt}
Let $\Omega$, $H_0$, $V$ and $\mathcal{F}$ satisfy Hypotheses \ref{cases}.
Then, there exists a (countably additive, complex)
measure $\nu_1=\nu_{1,H_0,V}$ on $\Omega$ such that
\begin{align*}
\|\nu_1\|\leq \min\left\{\tau_\Ical\big(|\re(V)|\big)+\tau_\Ical\big(|\im(V)|\big),\,\|V\|_\Ical\right\}
\end{align*} and
\begin{align*}
\tau_\Ical\big(f(H_0+V)-f(H_0)\big)=\int_\Omega f'(\la)\,\nu_1(d\la),
\end{align*}
for all $f\in\mathcal{F}$.
If Hypotheses \ref{cases}\eqref{i} are satisfied, then the measure $\nu_1$ is real and unique.
\end{thm}

When $\Ical=S^1$, the measure $\nu_1$ is absolutely continuous, but when $\Ical$ is the dual Macaev ideal (with the Dixmier trace), the measure $\nu_1$ can be of any type \cite[Theorem 4.4]{dsD}. Moreover, we do not have an explicit formula for $\nu_1$ in case of a general trace $\tau_\Ical$. 
Derivation of an explicit formula for $\nu_1$ in case $\Ical=S^1$, $H_0=H_0^*$, and $V=V^*$ relies on the fact that $\Tr\big(E_{H_0}(\cdot)V\big)$ is a (countably-additive) measure,
while the set function $\Tr_\omega\big(E_{H_0}(\cdot)V\big)$ can fail to be countably-additive (see \cite[Section 3]{dsD}).

As another consequence of singularity of $\Tr_\omega$ (and, more generally, of every trace satisfying $\tau_\Ical(\Ical^2)=\{0\}$), we have the following linearization formula. 

\begin{thm}
Assume Hypotheses \ref{cases} and assume $\tau_\Ical(\Ical^2)=\{0\}$. Then,
\begin{equation*}
\tau_\Ical\big(f(H_0+V)-f(H_0)\big)=\tau_\Ical\big(f'(H_0)V\big).
\end{equation*}
\end{thm}

Below we consider perturbations in the ideal $\Ical^{1/2}=\big\{A\in\Mcal:\;|A|^2\in\Ical\big\}$ and impose an additional natural assumption $\|AB\|_\Ical\leq\|A\|_{\Ical^{1/2}}\|B\|_{\Ical^{1/2}}$, which, in particular, holds for the ideals $S^1$, $\Lcal^1$, and $\Lcal^{(1,\infty)}$. 

\begin{hyp}
\label{2cases}
Consider a set $\Omega$, a closed, densely defined operator $H_0$ affiliated with $\Mcal$,
$V\in\Ical^{1/2}$ and a set $\mathcal{F}$ of functions
that satisfy one of the following assertions:
\begin{enumerate}[(i)]
\item \label{2ii} $\Omega=\Reals$, $\im(H_0)\geq 0$, $\im(H_0+V)\geq 0$, and
\[
\mathcal{F}=\text{\rm span}\left\{\la\mapsto (z-\la)^{-k}:\, k\in\Nats,\, \im(z)<0\right\};
\]
\item \label{2iii} $\Omega=\mathbb{T}$, $\|H_0\|\leq 1$, $\|H_0+V\|\leq 1$, and
$\mathcal{F}$ is the set of all functions that are
analytic on discs centered at $0$ and of radius strictly larger than $1$.
\end{enumerate}
\end{hyp}

\begin{thm}
\label{2mt}
Let $\Omega$, $H_0$, $V$ and $\mathcal{F}$ satisfy Hypotheses \ref{2cases}.
Then, there exists a (countably additive, complex)
measure $\nu_2=\nu_{2,H_0,V}$ on $\Omega$ such that
\begin{align*}
\|\nu_2\|\leq \frac12\,\tau_\Ical(|V|^2)
\end{align*} and
\begin{align*}
\tau_\Ical\left(f(H_0+V)-f(H_0)-\frac{d}{dt}\bigg|_{t=0}f(H_0+tV)\right)=\int_\Omega f''(\la)\,\nu_2(d\la),
\end{align*}
for every $f\in\mathcal{F}$.
\end{thm}

\begin{thm}
\label{2imt}
Suppose $\tau_\Ical(\Ical^{3/2})=\{0\}$.
Either assume Hypotheses \ref{2cases} or else take $H_0=H_0^*\in\Mcal$, $V=V^*\in\Ical^{1/2}$, and
$\mathcal{F}=C^4(\Reals)$.
Then, for every $f\in\mathcal{F}$,
\begin{align*}
\tau_\Ical\left(f(H_0+V)-f(H_0)-\frac{d}{dt}\bigg|_{t=0}f(H_0+tV)\right)=
\frac12\,\tau_\Ical\left(\frac{d^2}{dt^2}\bigg|_{t=0}f(H_0+tV)\right).
\end{align*}
\end{thm}

The major components in the proofs of Theorems \ref{mt} and \ref{2mt} are analogs of the estimates
\eqref{dest}, which hold due to the continuity of $\tau_\Ical$ with respect to $\|\cdot\|_\Ical$. However, presence of a singular component in the trace $\tau_\Ical$ requires more careful treatment of the operator derivatives than in the case of the normal trace $\Tr$.

\bibliographystyle{plain}

\end{document}